\documentclass{article}
\usepackage[tbtags]{amsmath}
\usepackage{amsfonts,amssymb,mathrsfs,amscd,comment}

\newcommand{\beq}{\begin{equation}}
\newcommand{\be}{\begin{equation}}
\newcommand{\ee}{\end{equation}}
\newcommand{\ba}{\begin{eqnarray}}
\newcommand{\ea}{\end{eqnarray}}

\newcommand{\lab}[1]{\label{#1}}

\newcommand{\CC}{\mathbb C}

\newcommand{\Z}{\mathbb Z}

\begin{document}

\title{\bf\large
Askey-Wilson continued fraction at $q^N=1$ }%

\author{\bf \normalsize A. R. Ismailov$^{1}$ and V. P. Spiridonov$^{1,2}$}

\date{}

\maketitle

\makeatletter
\renewcommand{\@makefnmark}{}
\makeatother
\footnotetext{
$\makebox[-2.1em]{}$
This study has been partially funded within the framework of the HSE University Basic
Research Program.\\
$^{1}$ NRU HSE, Moscow; ${}^2$ BLTP JINR, Dubna}

Consider the three-term recurrence relation
\begin{equation}
X_{n+1}(x)+ v_n X_n(x) + u_n X_{n-1}(x)=xX_n(x), \quad n\in\mathbb{Z}_{> 0},\; x\in\mathbb{C},
\label{ttrAWP}\end{equation}
which under the initial conditions  $X_0(x)=1,\, X_1(x)=x-v_0$ and the choice
\begin{eqnarray}\lab{uncn} &&
u_n=\xi_n \eta_{n-1}, \quad v_n= a +a^{-1}-\xi_n -\eta_n,
\\ \nonumber &&
\xi_n= \frac{a (1-q^n)(1-bcq^{n-1})(1-bdq^{n-1})(1-cdq^{n-1})}
{(1-gq^{2n-2})(1-gq^{2n-1})}, \quad g=abcd,
\\ \nonumber &&
\eta_n= \frac{a^{-1}(1-abq^n)(1-acq^n)(1-adq^n)(1-gq^{n-1})}
{(1-gq^{2n})(1-gq^{2n-1})},\quad a,b,c,d,q\in\CC,
\end{eqnarray}
generates the monic Askey-Wilson polynomials \cite{AW}, $X_n:=P_n(x),\, n\geq 0$.
These are the most general classical orthogonal polynomials which have found
important applications in algebraic combinatorics and mathematical physics.
Explicitly they are expressed  in terms of the $q$-hypergeometric series ${_4}\varphi_3$,
\begin{equation}
P_n(x)= \frac{a^{-n}(ab;q)_n(ac;q)_n(ad;q)_n}{(gq^{n-1};q)_n}
\;{_4}\varphi_3 \left({q^{-n},gq^{n-1},at,a/t
\atop ab, ac, ad} ;q, q \right),
\label{Pn} \end{equation}
where $x=t+t^{-1}$ and $(a;q)_n=\prod_{m=0}^{n-1}(1-aq^m)$ denotes the $q$-Pochhammer symbol.

Initial conditions $X_0^{(1)}=0,\, X_1^{(1)}=1$ in \eqref{ttrAWP} produce the associated
Askey-Wilson polynomials, $X_{n}^{(1)}:=P_{n-1}^{(1)}(x),\, n\ge 1$. It is known that the ratio
$P_{n-1}^{(1)}/P_{n}$ is equal to the continued fraction \cite{JT},
\begin{equation}
F_{n}(x)=\frac{P_{n-1}^{(1)}(x)}{P_{n}(x)}={1\over\displaystyle x-v_0 - {\strut u_1
\over\displaystyle x-v_1
- \dots {\atop \displaystyle -\frac{u_{n-1}}{x-v_{n-1}} }}},
\quad n\in\Z_{> 0}.
\label{cf-chain} \end{equation}
If the coefficient  $u_N=0$ for some $N>0$, then the continued fraction $F_n(x)$,
for all $n\ge N$, is truncated at the $N$-th place and, therefore, coincides with $F_N(x)$.
The main result of the present work consists in the computation of the explicit form of
$F_{N}(x)$, when $u_N=0$ due to the choice of $q$ equal to a root of unity, $q^N=1$.
Let us denote
$$
E_N={a^N+b^N+c^N+d^N -(abc)^N-(bcd)^N-(abd)^N-(acd)^N \over 1-(abcd)^N},
$$
and $r$ --- one of the solutions of the equation $r^N + r^{-N}= E_N$.

{\bf Theorem.}
{\em Let $q=e^{2\pi i M/N},\, (M,N)=1$, and the parameters $a,b,c,d$ satisfy the cons\-traints
$E_N \ne \pm 2$ and  $g, ab, ac, ad,$ $bc,$ $bd, cd  \ne q^k, \, k\in\mathbb{Z}$.
Then the rational function  $F_N(x)$ has the following representation:
\begin{equation}
F_N(x)=
\frac{\kappa}{x-r-1/r}\frac{\displaystyle {_9}\varphi_8\left({q^{-(N-1)},rt,r/t,qr,-qr,ar,br,cr,dr
\atop qr/t,qrt,r,-r,qr/a,qr/b,qr/c,qr/d} ;q, \frac{q^2}{g} \right)}
{\displaystyle {_4}\varphi_3 \left({q^{-(N-1)},gq^{-2},ar,a/r \atop ab, ac, ad} ;q,q \right)},
\label{FNfin}\end{equation}
where
$$
\kappa=a^{N-1}\frac{(1-(bc)^N)(1-(cd)^N)(1-(bd)^N)(1-g/q)
(1-g/q^2)(r-r^{-1})}{(1-g^N)^2 (1-bc/q)(1-cd/q)(1-bd/q)(r^N-r^{-N})}.
$$
}

{\bf Proof.} The indicated choice of $q$ and the conditions $g, ab, ac, ad,$ $bc,$
$bd, cd  \ne q^k, \, k\in\mathbb{Z}$, guarantee that the coefficients of the recurrence relation
are finite, $u_n\neq 0$ for $0<n<N$, and $u_N=0$. In this case poles of the rational function
$F_N(x)$ are determined by the zeros of the polynomial $P_N(x)=t^N + t^{-N} - E_N$,
the explicit form of which is fixed by the limiting transition $\mu\to 1^-$ for $q=\mu e^{2\pi i M/N},\,
(M|N)=1,\, \mu\in\mathbb{R}$. The equality $P_N(x)=0$ holds when $t^N=E_N/2 \pm \sqrt{E_N^2/4 -1}$,
i.e., the roots of this polynomial have the form $x_k= rq^k+r^{-1}q^{-k}, \, k=0,1,\dots,N-1$, where $r$
is one of the admissible values of $t$. Correspondingly, $x_k\neq x_l,\, k\neq l$,
if $E_N\neq\pm 2$, which is assumed in the following.
The condition $P_N(x_k)=0$ generates a finite system of Askey-Wilson polynomials which was investigated
in \cite{SZroots}, and we use some of the results derived in that work.

Because zeros of $P_N(x)$ are simple, one has the partial fraction decomposition
\begin{equation}
F_N(x)=\sum_{k=0}^{N-1} \frac{w_k}{x-x_k}, \qquad
w_k = \frac{P^{(1)}_{N-1}(x_k)}{P_{N}'(x_k)},
\label{F_frac} \end{equation}
where the coefficients $P_N'(x_k)$ have the following form
\begin{equation}
P_N^{\prime}(x_k)= N\frac{r^N-r^{-N}}{rq^k-r^{-1}q^{-k}}.
\lab{P'}\end{equation}
The Casorati determinant of the chosen solutions of the recurrence relation
\eqref{ttrAWP} $P_n(x)$ and $P_n^{(1)}(x)$ is exactly computable
\begin{equation}
P_{n-1}(x) P^{(1)}_{n-1}(x) - P_n(x)P^{(1)}_{n-2}(x) = h_{n-1}, \quad h_n=\prod_{k=1}^n u_k, \; h_0=1.
\label{Cas} \end{equation}
For $n=N$ and $x=x_k$ one has $P_N(x_k)=0$, so that \eqref{Cas} leads to the expression
\begin{equation}
w_k = \frac{h_{N-1}}{P_{N}'(x_k) P_{N-1}(x_k)}.
\label{wk1} \end{equation}

According to \cite{SZroots}, the coefficients $w_k$ are the weights in the discrete orthogonality relation
$\sum_{k=0}^{N-1} P_n(x_k)P_m(x_k)w_k = h_n \delta_{nm}, \, n,m=0,1,\ldots, N-1.$
The weights $w_k$ were computed in \cite{SZroots} by a different method
\be
w_k= w_0 \frac{(1-r^2q^{2k})(ar;q)_k(br;q)_k(cr;q)_k(dr; q)_k}
{(1-r^2)(qr/a;q)_k(qr/b;q)_k(qr/c;q)_k(qr/d; q)_k}\Big(\frac{q}{g}\Big)^k,
\lab{wk2} \end{equation}
where $w_0$ is found from the previous expression \eqref{wk1}. Substituting the representation
\begin{equation}
x-x_k =\frac{t+1/t-r-1/r}{q^k}\frac{(rqt;q)_k(rq/t;q)_k}{(rt;q)_k(r/t;q)_k}
\label{x}\end{equation}
into the expansion \eqref{F_frac} and transforming the emerging $q$-hypergeometric
series to a standard form we obtain the needed expression for the Askey-Wilson continued fraction
\eqref{FNfin}. This is a new explicit continued fraction which is expressed as a ratio of two
terminating  $q$-hypergeometric series.


\begin{thebibliography}{000000}

\bibitem{AW} R.~ Askey and J.~ Wilson, Mem. Amer. Math. Soc.
{\bf 54} (1985), 1-55.

\bibitem{JT} W.~ B.~ Jones, W.~ J.~ Thron, {\em Continued Fractions: Analytic
Theory and Applications}, Addison-Wesley, 1980.

\bibitem{SZroots} V.~ Spiridonov, A.~ Zhedanov, Duke Math. J. {\bf 89} (1997), 283-305.

\end{thebibliography}
\end{document}